\documentclass[12pt]{amsart}
\usepackage[latin1]{inputenc}
\usepackage{indentfirst}
\usepackage[active]{srcltx}
\usepackage{amsfonts}
\usepackage{graphicx}
\usepackage{amsmath}
\usepackage{amssymb}
\usepackage{latexsym}
\usepackage{amscd}
\usepackage[all]{xy}
\usepackage{xcolor}
\usepackage{apptools,chngcntr} 
\newcommand\mut[1]{\ignorespaces}
\usepackage{verbatim}
\usepackage{url}

\newcommand{\lub}{\mbox{lub}}

\newcommand{\chav}[1]{\ensuremath{\left\{#1\right\}}}

\newtheorem{theorem}{Theorem}[section]
\newtheorem{definition}[theorem]{Definition}
\newtheorem{lemma}[theorem]{Lemma}
\newtheorem{corollary}[theorem]{Corollary}
\newtheorem{proposition}[theorem]{Proposition}
\newtheorem{remark}[theorem]{Remark}
\newtheorem{example}[theorem]{Example}

\allowdisplaybreaks
\begin{document}

\title[]{Weierstrass semigroup at $m+1$ rational points in maximal curves which cannot be covered by the Hermitian curve}

\author{Alonso Sepúlveda Castellanos and Maria Bras-Amorós}

\maketitle

\begin{abstract}
We determine the Weierstrass semigroup $H(P_\infty,P_1,\ldots,P_m)$ at se\-ve\-ral rational points on the maximal curves which cannot be covered by the Hermi\-tian curve introduced in \cite{TTT}. Furthermore, we present some conditions to find pure gaps. We use this semigroup to obtain AG codes with better relative parameters than comparable one-point AG codes arising from these curves.
\end{abstract}


\section{Introduction}
Let $\mathcal{X}$ be a non-singular, projective, irreducible, algebraic curve of genus $g\geq 1$ over a finite field $\mathbb{F}_{q^{2}}$ with genus
$g\left( \mathcal{X}\right) $ and $\#\mathcal{X}\left( {\mathbb{F}_{q^{2}}}\right)$ rational points. The Hasse-Weil bound states that 
\[
\left\vert \#\mathcal{X}\left( {\mathbb{F}_{q^{2}}}\right)-\left( q^2+1\right) \right\vert \leq 2qg\left( \mathcal{X}\right)\,.
\]
In the case that $\#\mathcal{X}\left({\mathbb{F}_{q^{2}}}\right)=q^2+1+2qg\left( \mathcal{X}\right)$, the curve $\mathcal{X}$ is called an ${\mathbb{F}_{q^{2}}}$-{\it maximal curve}. Maximal curves have been widely studied \cite{FGT},\cite{KT}. We know that every curve that is covered by an ${\mathbb{F}_{q^{2}}}$-maximal curve also turns out to be an ${\mathbb{F}_{q^{2}}}$-maximal curve, see \cite{Lan}.

In \cite{GS2}, Garcia and Stichtenoth presented an example of a maximal curve that is not Galois, covered by the maximal Hermitian curve.  Later, Giulietti and Korchm\'{a}ros \cite{GK}  presented a family of maximal curves which cannot be covered by the Hermitian curve, the $GK$-curve. Garcia, Guneri and Stichtenoth \cite{GGS} present a generalization of this curve, the $GGS$-curves. In \cite{TTT}, Tafazolian, Teher\'{a}n and Torres presented two families of maximum curves that could not be covered by the Hermitian curve, the $\mathcal{X}_{a,b,n,s}$ and $\mathcal{Y}_{n,s}$ curves, and in \cite{BM} Beelen and Montanucci construct another generalization of the $GK$-curve.      

Curves with many rational points have been investigated to construct error-correcting codes, the so-called {\it algebraic geometric codes} (AG codes), introduced by Goppa \cite{goppa1}, \cite{goppa2}. An important part of the study of AG codes are the Weierstrass semigroups on rational points of the curve because there exists a close connection between the parameters of one-point AG codes and their duals with the Weierstrass semigroup at one point on the curve. Weierstrass semigroups have been used to analyze the minimum distance, as in 
\cite{HoLiPe,KiPe,FeRa:d}, 
to analyze the code redundancy and to determine improved codes, as in \cite{PeTo,BrOS:cc,Br:redundancy}, to bound the code length, as in \cite{GeMa}, or to analyze the weight hierarchy and the generalized Hamming weights, as in 
\cite{HePe,FaMu,BLV}.
This way, an effort was put to explicitly compute Weierstrass semigroups of particular families of curves \cite{CaFa,CaFaMu,PeStTo}. The case of one-point, two-point and multi-point AG codes on the GK maximal curves were studied in \cite{FG}, \cite{CT1} and \cite{CT},\cite{BMZ}, respectively. In \cite{BMZ1}, Bartoli, Montanucci and Zini examined one-point AG codes from the GGS curves, and in \cite{HuYang}, Hu and Yang studied multi-point AG codes on the GGS curves. They explicited bases for Riemann-Roch spaces using a related set of lattice points and considered the properties from GGS curves to characterize the Weierstrass semigroup and pure gaps in many rational places. In this contribution we determine the Weierstrass semigroup in many points of the $\mathcal{X}_{a,b,n,s}$ and $\mathcal{Y}_{n,s}$ curves, using the concept of {\em discrepancy} introduced by Duursma and Park in \cite{Duursma}. 

This paper is organized as follows. Section 2 and Section 3 contain general results on Weierstrass semigroups, basic facts related to AG codes and the basic properties of the $\mathbb{F}_{q^{2n}}$-maximal curves $\mathcal{X}_{a,b,n,s}$ and $\mathcal{Y}_{n,s}$. In Section 4 we describe some generators of the Weierstrass semigroup at rational points of those curves of the form $P_{(\alpha,\beta,0)}$ and in Appendix A and Appendix B we prove that the genera of the semigroups generated by these generators coincides with the genus of the curve, thus proving that the  described generators are, indeed, {\em all} the generators of the Weierstrass semigroup.
In Section 5 we determine the minimal generating set for the Weierstrass semigroup $H(P_\infty,P_1,\ldots,P_m)$, where $P_1,\ldots,P_m$ are rational points on the curves with $1\leq m\leq q$. Finally, in Section 6 we present some results about pure gaps and AG codes and we illustrate them with an example of a code with 
better relative parameters than comparable one-point AG codes, and an example of a quasi perfect code.

\section{Preliminaries}

Let $\mathcal{X}$ be a non-singular, projective, irreducible, algebraic curve of genus $g \geq 1$ over a finite field $\mathbb{F}_{q}$.
Fix $m$ distinct rational points $P_1,\ldots,P_m$ on $\mathcal{X}$ and let $H(P_1,\ldots,P_m)$ be the Weiertrass semigroup at  $P_1,\ldots,P_m$. Define a partial order $\preceq$ on $\mathbb{N}_0^m$ by $(n_1,\ldots,n_m)\preceq (p_1,\ldots,p_m)$ if and only if $n_i\leq p_i$ for all $i$, $1\leq i\leq m$.
For ${\bf u}_1,\ldots,{\bf u}_t\in \mathbb{N}_0^{m}$, where, for all $k$, ${\bf u}_k = (u_{k_{1}}, \ldots , u_{k_{m}})$, we define the \emph{least upper bound} ($\lub$) of the vectors ${\bf u}_1,\ldots,{\bf u}_t$ in the following way: \[\lub\{{\bf u}_1,\ldots,{\bf u}_t\}=(\max\{{ u_{1_1}},\ldots,{ u_{t_1}}\},\ldots, \max\{{ u_{1_m}},\ldots,{ u_{t_m}}\} )\in \mathbb{N}_0^{m}.\] For $\mathbf{n}=(n_1,\ldots,n_m)\in \mathbb{N}_0^{m}$ and $i \in \{ 1,\ldots , m\}$, we set
$$
\nabla_i (\mathbf{n}):=\{ (p_1, \ldots , p_m) \in H(P_{1}, \ldots, P_m) \mbox{ ; } p_i=n_i \}.
$$

\begin{proposition}\label{minimal}\cite[Proposition 3]{gretchen1}
Let $\mathbf{n}=(n_1,\ldots,n_m)\in \mathbb{N}_0^{m}$. Then $\mathbf{n}$ is minimal, with respect to $\preceq$, in $\nabla_i(\mathbf{n})$ for some $i$, $1 \leq i \leq m$, if and only if $\mathbf{n}$ is minimal in $\nabla_i(\mathbf{n})$ for all $i$, $1 \leq i \leq m$.
\end{proposition}


\begin{proposition}\label{lubH}\cite[Proposition 6]{gretchen1}
Suppose that $1 \leq t \leq m \leq q$ and ${\bf u}_1,\ldots,{\bf u}_t\in H(P_{1}, \ldots , P_{m})$. Then $\lub\{{\bf u}_1,\ldots,{\bf u}_t\} \in H(P_{1}, \ldots , P_{m})$.
\end{proposition}


\begin{definition}\label{defnabla}
Let $\Gamma(P_{1})=H(P_1)$ and, for $m\geq 2$, define
$$
\Gamma(P_{1}, \ldots, P_{m}):=\{{\bf n}\in \mathbb{N}^m: \mbox{ for some } i, 1\leq i\leq m, {\bf n} \mbox{ is minimal in } \nabla_i ({\bf n})\}.
$$
\end{definition}


\begin{lemma}\label{mPoints} \cite[Lemma 4]{gretchen1}
For $m \geq 2$, $\Gamma(P_{1}, \ldots, P_{m}) \subseteq G(P_1)\times\cdots\times G(P_m)\;,$ where $G(P_i)$ is the set of gaps in $P_i$.
\end{lemma}

In \cite{gretchen1}, Theorem 7, it is shown that, if $2\leq m \leq q$, then $H(P_1,\ldots,P_m)= $
\[ \left\{\begin{array}{cl} \lub\{{\bf u}_1,\ldots,{\bf u}_m\}\in \mathbb{N}_0^m: & {\bf u}_i \in \Gamma(P_{1}, \ldots, P_{m}) \\& \mbox{ or } (({\bf u}_i)_{j_1}, \dots, ({\bf u}_i)_{j_k}) \in \Gamma(P_{j_1}, \ldots , P_{j_k}) \\ & \mbox{ for some } \{j_1,\ldots,j_k\}\subset\{1,\ldots,m\}  \mbox{ such that } \\ & j_1<\cdots<j_k \mbox{ and } ({\bf u}_i)_{j_{k+1}} = \cdots = ({\bf u}_i)_{j_m}=0, \\ & \mbox{ where }  \{j_{k+1},\ldots,j_m\}=\{1,\ldots,m\} \setminus \{j_1,\ldots,j_k\} \end{array}\right\}\;.\]

Therefore, the Weierstrass semigroup $H(P_1,\ldots,P_m)$ is completely determined by the sets $\Gamma(P_1,\ldots,P_\ell)$ with $1\leq \ell\leq m$. In \cite{gretchen1}, Matthews called the set $\Gamma(P_1,\ldots,P_m)$ a \emph{minimal generating set} of $H(P_1,\ldots,P_m)$.

\section{The curves $\mathcal{X}_{a,b,n,s}$ and $\mathcal{Y}_{n,s}$}\label{scurve1}

\subsection{The curve $\mathcal{X}_{a,b,n,s}$}\label{curve1}

Let $a,b,s\geq 1,n\geq 3$ be integers such that $n$ is odd. Let $q=p^a$ be a power of a prime number $p$, $b$ is a divisor of $a$, $s$ is a divisor of $(q^n+1)/(q+1)$ and $c\in \mathbb{F}_{q^2}$ with $c^{q-1}=-1$. We define the curve $\mathcal{X}_{a,b,n,s}$ over $\mathbb{F}_{q^{2n}}$ by the affine equations 
\begin{equation}
cy^{q+1}=t(x):=\sum_{i=0}^{a/b-1}x^{p^{ib}}\;\;\mbox{and}\;\; z^{M}=y^{q^2}-y\;,
\end{equation}
where $M=\dfrac{q^n+1}{s(q+1)}$. This curve is absolutely irreducible, nonsingular, and maximal over $\mathbb{F}_{q^{2n}}$ of genus $g=\dfrac{q^{n+2}-p^bq^n-sq^3+q^2+(s-1)p^b}{2sp^b}$. From Theorem 3.5 in \cite{TTT}, the curve $\mathcal{X}_{a,b,n,1}$ cannot be Galois covered by the Hermitian curve $\mathcal{H}_n: v^{q^n+1}=u^{q^n}+u$ over $\mathbb{F}_{q^{2n}}$ provided that $b<a$. A plane model of $\mathcal{X}_{a,b,n,s}$ is given by the equation \[cz^{\frac{q^n+1}{s}}=t(x)(t(x)^{q-1}+1)^{q+1}\;.\]

Let $\mathcal{X}_{a,b,n,s}(\mathbb{F}_{q^{2n}})$ be the set of $\mathbb{F}_{q^{2n}}$-rational points of $\mathcal{X}_{a,b,n,s}$, and we will denote a rational point $P=(\alpha,\beta,\gamma)\in\mathcal{X}_{a,b,n,s}(\mathbb{F}_{q^{2n}})$ by $P_{(\alpha,\beta,\gamma)}$, whereas $P_0=(0,0,0)$. Let $P_\infty$ be the unique common pole of the functions $x,y,z$ which define the function field of $\mathcal{X}_{a,b,n,s}$, then we have the following divisors:
\begin{equation*}
(x-\alpha) = (q+1)MP_{(\alpha,\beta,0)}-(q+1)MP_\infty\;, \mbox{with }t(\alpha)=c\beta^{q+1} \mbox{ and }\beta\in\mathbb{F}_{q^2}\;,
\end{equation*}
\begin{equation}\label{divisorscurve1}
(y-\beta) = \sum_{i=1}^{q/p^b}MP_{(\alpha_i,\beta,0)}-\dfrac{q}{p^b}MP_\infty\;,\mbox{with }t(\alpha_i)=c\beta^{q+1} \mbox{ and }\beta\in\mathbb{F}_{q^2}\;,
\end{equation}
\begin{equation*}
(z)  =  \sum_{j=1}^{q^2}\sum_{i=1}^{q/p^b}P_{(\alpha_i,\beta_j,0)}-\dfrac{q^3}{p^b}P_\infty\;, \mbox{with }\beta_j\in \mathbb{F}_{q^2}\mbox{ and }c\beta_j^{q+1}=t(\alpha_i), \forall i=1,\ldots,q/p^b\;.
\end{equation*}
From \cite[Proposition 5.1]{TTT}, we have that $H(P_\infty)=\langle \frac{q}{p^b}M, \frac{q^3}{p^b}, (q+1)M \rangle$. 

\subsection{The curve $\mathcal{Y}_{n,s}$}\label{curve2}

Let $n\geq 3$ be an odd integer, let $q$ be a prime power, and let $s\geq 1$ be a divisor of $\frac{q^n+1}{q+1}$. We define the curve $\mathcal{Y}_{n,s}$ over $\mathbb{F}_{q^{2n}}$ by the affine equations 
\begin{equation}
y^{q+1}=x^q+x\;\mbox{ and }\; z^{M}=y^{q^2}-y\;,
\end{equation}
where $M=\frac{q^n+1}{s(q+1)}$. This curve is maximal over $\mathbb{F}_{q^{2n}}$ of genus $g(\mathcal{Y}_{n,s})=\frac{q^{n+2}-q^n-sq^3+q^2+s-1}{2s}$. From \cite[Theorem 4.4]{TTT}, we know that the curve $\mathcal{Y}_{3,s}$ cannot be covered by the Hermitian curve $H_3$ over $\mathbb{F}_{q^6}$, in case $q>s/(s+1)$. A plane model of $\mathcal{Y}_{n,s}$ is given by the equation \[z^{\frac{q^n+1}{s}}=(x^q+x)((x^q+x)^{q-1}-1)^{q+1}\;.\]

Let $\mathcal{Y}_{n,s}(\mathbb{F}_{q^{2n}})$ be the set of $\mathbb{F}_{q^{2n}}$-rational points of $\mathcal{Y}_{n,s}$, and we will denote a rational point $P=(\alpha,\beta,\gamma)\in\mathcal{Y}_{n,s}(\mathbb{F}_{q^{2n}})$ by $P_{(\alpha,\beta,\gamma)}$, whereas $P_0=(0,0,0)$. Let $P_\infty$ be the unique common pole of the functions $x,y,z$ which define the function field of $\mathcal{Y}_{n,s}$, then we have the following divisors:
\[(x-\alpha) = (q+1)MP_{(\alpha,\beta,0)}-(q+1)MP_\infty\;, \mbox{with } \alpha^q+\alpha=\beta^{q+1} \mbox{ and }\beta\in\mathbb{F}_{q^2}\;,\]
\begin{equation}\label{divisorscurve2}
(y-\beta) =  \sum_{i=1}^{q}MP_{(\alpha_i,\beta,0)}-qMP_\infty\;,\mbox{with }\alpha_i^q+\alpha_i=\beta^{q+1} \mbox{ and }\beta\in\mathbb{F}_{q^2}\;,
\end{equation}
\begin{equation*}
(z)  =  \sum_{j=1}^{q^2}\sum_{i=1}^{q}P_{(\alpha_i,\beta_j,0)}-q^3P_\infty\;, \mbox{with }\beta_j\in \mathbb{F}_{q^2}\mbox{ and }\beta_j^{q+1}=\alpha_i^q+\alpha_i, \forall i=1,\ldots,q\;.
\end{equation*}
From \cite[Proposition 5.1]{TTT}, we have that $H(P_\infty)=\langle qM, q^3, (q+1)M \rangle$. 
For $s=1$, we have that $\mathcal{Y}_{n,1}=GGS(\mathcal{X})$-curves.

\section{The Weierstrass Semigroup $H(P_{(\alpha,\beta,0)})$}

\begin{proposition} Let $b<a$.

The Weierstrass semigroup at $P_{(\alpha,\beta,0)}\in\mathcal{X}_{a,b,n,1}(\mathbb{F}_{q^{2n}})$  is
\[H(P_{(\alpha,\beta,0)})=\left\langle q^n+1-iM-j: 0\leq i\leq p^b, 0\leq j\leq q^{n-3}p^b-i\dfrac{q^{n-2}+1}{q+1}\right\rangle\]

The Weierstrass semigroup at $P_{(\alpha,\beta,0)}\in\mathcal{Y}_{n,1}(\mathbb{F}_{q^{2n}})$ is
\[H(P_{(\alpha,\beta,0)})=\left\langle q^n+1-iM-j: 0\leq i\leq 1, 0\leq j\leq q^{n-3}-i\dfrac{q^{n-2}+1}{q+1}\right\rangle\]
\end{proposition}

\begin{proof}
For $n=3$, we have that $\left( \dfrac{(y-\beta)^iz^j}{x-\alpha} \right)_\infty=(q^3+1-iM-j)P_{(\alpha,\beta,0)}\;,$ for $0\leq i+j\leq p^b$. Then the numerical semigroup containing all linear combinations with non-negative integer coefficients of these values, i.e, $S=\langle q^3+1-iM-j: 0\leq i+j\leq p^b \rangle$ is contained in $H(P_{(\alpha,\beta,0)})$. By the proof in Appendix~\ref{a:S}, we have that the genus of the semigroup $S$ is equal to $g(\mathcal{X})$ and therefore the assertion follows.

For $n\geq 5$ odd, we have that $\left( \dfrac{(y-\beta)^iz^j}{x-\alpha} \right)_\infty=(q^n+1-iM-j)P_{(\alpha,\beta,0)}\;,$ for $0\leq iM+jq^2\leq q^{n-1}p^b$. Now, the numerical semigroup  $S'=\langle q^n+1-iM-j: 0\leq iM+jq^2\leq q^{n-1}p^b\rangle\subseteq H(P_{(\alpha,\beta,0)})$.
By the proofs in Appendix~\ref{a:Sprime}, we have that the genus of the semigroup $S'$ is equal to $g(\mathcal{X})$ and therefore the assertion follows. We can observe that if $j=0$ then the maximal value for $i$ is $p^b$ in $0\leq iM+jq^2\leq q^{n-1}p^b$. Furthermore, \[ j\leq \dfrac{q^{n-1}p^b-i(q^n+1)/(q+1)}{q^2}\leq q^{n-3}p^b-i\dfrac{q^n+1}{q^2(q+1)}\leq q^{n-3}p^b-i\dfrac{q^{n-2}+1}{q+1}+i\dfrac{q-1}{q^2}\;.\] The range of the parameters $i,j$ follows now from the inequality $i(q-1)/q^2<1$.
\end{proof}

\begin{remark}
In \cite[Proposition 4.3]{BMZ1}, Bartoli, Montanucci and Zini calculate the Weierstrass semigroup $H(P_{(\alpha,\beta,0)})$ for the curves $\mathcal{Y}_{n,1}$ in a different way. They observed that this semigroup is independent of the choice of $\alpha$ and $\beta$ by \cite[Lemma 8.1]{BMZ1}.
\end{remark}

\begin{example}\label{ex1curve1} For $s=1,n=3,p=2,a=2,b=1$, and $c=1$, we have that $q=4, M=13$, and the affine equations of the curve $\mathcal{X}_{2,1,3,1}$ are $y^5=x+x^2$ and $z^{13}=y^{16}-y$ with genus $g=212$. In this case $H(P_\infty)=\langle 65,32,26 \rangle$ and $H(P_{(\alpha,\beta,0)})=\langle 65,64,63,52,51,39 \rangle$.
\end{example}

\begin{example}\label{ex2curve1}
For $s=1,n=5,p=2,a=2,b=1$ and $c=1$, we have that $q=4$ and $M=205$, and the affine equations of the curve $\mathcal{X}_{2,1,5,1}$ are $y^5=x+x^2$ and $z^{205}=y^{16}-y$ with genus $g=3572$. In this case $H(P_\infty)=\langle 1025,410,32 \rangle$ and $H(P_{(\alpha,\beta,0)})=\langle 1025,\ldots,993,820,\ldots,801,615,\ldots,609 \rangle$.
\end{example}

\begin{example}\label{ex1curve2}
For $s=1,n=3,q=3$, we have $M=7$ and the affine equations of the curve $\mathcal{Y}_{3,1}$ are $y^4=x^3+x$ and $z^{7}=y^{9}-y$ with genus $g=99$. In this case $H(P_\infty)=H(P_{(\alpha,\beta,0)})=\langle 21,27,28 \rangle$.
\end{example}

\begin{example}\label{ex2curve2}
For $s=1,n=5,q=2$, we have $M=11$ and the affine equations of the curve $\mathcal{Y}_{5,1}$ are $y^3=x+x^2$ and $z^{11}=y^{4}-y$ with genus $g=46$. In this case $H(P_\infty)=\langle 33,22,8 \rangle$ and $H(P_{(\alpha,\beta,0)})=\langle 33,32,31,30,29,22,21 \rangle$.
\end{example}

\section{The Weierstrass semigroup at certain $m+1$ points on the curve $\mathcal{C}$}\label{Weierstrass-S}

Let $a,b,n,s=1,p,q=p^a,M=\dfrac{q^n+1}{q+1}$ be as above, and let us fix the following notation:

\begin{itemize}
	\item $\mathcal{C}$ denotes either the curve $\mathcal{X}_{a,b,n,1}$ in subsection \ref{curve1}, or the curve $\mathcal{Y}_{n,1}$ in subsection \ref{curve2}.
	
	\item $P_\infty\in \mathcal{C}(\mathbb{F}_{q^{2n}})$ is the unique common pole of the functions $x,y,z$ which define the function field of $\mathcal{C}$.
	
	\item $P_i:=P_{(\alpha_i,0,0)}\in \mathcal{C}(\mathbb{F}_{q^{2n}})$ for $i=1,\ldots,q/p^b$, and for $j=1,\ldots,(q^3-q)/p^b$, let  $Q_j=P_{(\alpha_j,\beta_j,0)}\in \mathcal{C}(\mathbb{F}_{q^{2n}})$ such that $\beta_j\neq 0$. 
\end{itemize}

In this section, we determine the Weierstrass semigroup $H(P_\infty,P_1,\ldots,P_m)$ for $1\leq m\leq q/p^b$ ($b=0$ when $\mathcal{C}=\mathcal{Y}_{n,s}$).

By the divisor of rational functions $(x-\alpha_\ell), y$ and $z$ given by \eqref{divisorscurve1} and \eqref{divisorscurve2}, we have the following equivalences:
\begin{equation}\label{eq1}
(q^n+1)P_\ell \sim (q^n+1)P_\infty\,,
\end{equation}
\begin{equation}\label{eq2}
MP_1+\cdots +MP_{q/p^b}\sim (q/p^b)MP_\infty\,,
\end{equation}
\begin{equation}\label{eq3}
P_1+\cdots+P_{q/p^b}+Q_1+\cdots+Q_{(q^3-q)/p^b}\sim (q^3/p^b)P_\infty\,.
\end{equation}

Let $1\leq m\leq q/p^b$, and let $1\leq k\leq M, 0\leq i\leq q$ and $j_\ell\geq 0$ be integers such that 

\[ \left( q^2-mp^b-p^b\sum_{\ell=1}^m j_\ell\right)(q^n+1)-iqM-kq^3 >0\,. \]

So, the divisor \[A'=\dfrac{1}{p^b}((q^n+1)(q^2-mp^b)-iqM-kq^3)P_\infty+\sum_{\ell=1}^m(iM+k)P_\ell\;\]
is effective and by \eqref{eq1}, we have that 

\begin{equation}\label{effectivedivisor}
\dfrac{1}{p^b}\left( \left( q^2-mp^b-p^b\sum_{\ell=1}^m j_\ell\right)(q^n+1)-iqM-kq^3 \right)P_\infty + \sum_{\ell=1}^m(j_\ell(q^n+1)+iM+k)P_\ell\sim A'\;.
\end{equation}

Next we state Duursma and Park's definition of discrepancy \cite[Section 5]{Duursma}.

\begin{definition}
A divisor $A'\in \mbox{Div}(\mathcal{X})$ is called a discrepancy for
two rational points $P$ and $Q$ on $\mathcal{X}$ if
$\mathcal{L}(A')\neq \mathcal{L}(A'-P)=\mathcal{L}(A'-P-Q)$ and
$\mathcal{L}(A')\neq \mathcal{L}(A'-Q)=\mathcal{L}(A'-P-Q)$.
\end{definition}

\begin{lemma}\cite[Lemma 2.6]{CT}\label{LemaCT} Let $\mathbf{n}=(n_1,\ldots,n_m)\in H(P_1,\ldots,P_m)$. Then $\mathbf{n}\in \Gamma(P_1,\ldots,P_m)$ if only if the divisor $A'=n_1P_1+\cdots+n_mP_m$ is a discrepancy with respect to $P$ and $Q$ for any two rational points $P,Q\in \{P_1,\ldots,P_m\}$.
\end{lemma}

\begin{lemma}\cite[Noether's Reduction Lemma]{Fulton} Let $D$ be a divisor, $P\in \mathcal{C}$ and let $K$ be a canonical divisor. If $\dim(\mathcal{L}(D))>0$ and $\dim(\mathcal{L}(K-D-P))\neq \dim(\mathcal{L}(K-D))$, then $\dim(\mathcal{L}(D+P))=\dim(\mathcal{L}(D))$.
\end{lemma}

\begin{proposition}\label{Discrepancy}
The divisor $A'$ is a discrepancy with respect to $P$ and $Q$ for any two distinct rational points $P,Q\in \{P_\infty,P_1,\ldots,P_m\}$.
\end{proposition}
\begin{proof}
From the equivalence in \eqref{effectivedivisor}, we have that there exists a function $f\in \mathcal{L}(A')$ with pole divisor equal to $A'$. Thus, $\mathcal{L}(A')\neq \mathcal{L}(A'-P)$ for every rational point $P\in \{P_\infty,P_1,\ldots,P_m\}$.

Now, we prove that $\mathcal{L}(K-A'+P)\neq \mathcal{L}(K-A'+P+Q)$, where $K$ is a canonical divisor. Let $K=\frac{1}{p^b}((q^n+1)(q^2-p^b)-q^3-p^b)P_\infty$, so
\[K+P+Q-A'=\frac{1}{p^b}((q^n+1)(q^2-p^b)-q^3-p^b)P_\infty+P+Q-\dfrac{1}{p^b}((q^n+1)(q^2-mp^b)-iqM-kq^3)P_\infty\] \[-\sum_{\ell=1}^m(iM+k)P_\ell=\frac{1}{p^b}((q^n+1)p^b(m-1)+iqM+(k-1)q^3-p^b)P_\infty+P+Q-\sum_{\ell=1}^m(iM+k)P_\ell\;.\]
Without loss of generality, we can assume that $P=P_\infty$ and $Q=P_1$. Thus, 
\[K+P_\infty+P_1-A'=\frac{1}{p^b}((q^n+1)p^b(m-1)+iqM+(k-1)q^3)P_\infty-(iM+k-1)P_1-\sum_{\ell=2}^m(iM+k)P_\ell\;,\] and we have that
\[z^{k-1}y^i(x-\alpha_2)\cdots (x-\alpha_m)\in \mathcal{L}(K+P_\infty+P_1-A')\setminus\mathcal{L}(K+P_1-A')\;.\]
So, $\mathcal{L}(A'-P_1)=\mathcal{L}(A'-P_\infty-P_1)$. Since $\mathcal{L}(A')\neq \mathcal{L}(A'-P_1)=\mathcal{L}(A'-P_\infty-P_1)$, and $\mathcal{L}(A')\neq \mathcal{L}(A'-P_\infty)$, it follows that $\mathcal{L}(A'-P_\infty)=\mathcal{L}(A'-P_\infty-P_1)$.

Now, if $P\neq P_\infty$ and $Q\neq P_\infty$, then we can suppose that $P=P_1$ and $Q=P_2$. In this case, we have that 
\[z^{k-1}y^i(x-\alpha_3)\cdots(x-\alpha_m)\in \mathcal{L}(K+P_1+P_2-A')\setminus\mathcal{L}(K+P_2-A')\,.\]

As above, we have that $\mathcal{L}(A'-P_2)=\mathcal{L}(A'-P_1-P_2)$ and that $\mathcal{L}(A'-P_1)=\mathcal{L}(A'-P_1-P_2)$. Therefore, the divisor $A'$ is a discrepancy with respect to $P$ and $Q$ for any two distinct rational points $P,Q\in \{P_\infty,P_1,\ldots,P_m\}$. 
\end{proof}

As a consequence of \eqref{effectivedivisor} and Proposition \ref{Discrepancy}, we have that the effective divisor $$A=\dfrac{1}{p^b}\left( \left( q^2-mp^b-p^b\sum_{\ell=1}^m j_\ell\right)(q^n+1)-iqM-kq^3 \right)P_\infty + \sum_{\ell=1}^m(j_\ell(q^n+1)+iM+k)P_\ell\,,$$ is also a discrepancy with respect to $P$ and $Q$ for any two distinct rational points $P,Q\in \{P_\infty,P_1,\ldots,P_m\}$.

\begin{theorem}\label{TeoGamma}
Let $a,b,n,s,p,q,M,P_\infty,P_1,\ldots,P_m$ be as above. For $1\leq m\leq q/p^b$, let 

\[\Gamma_{m+1}=\left\{ \left( \dfrac{1}{p^b}\left( \left( q^2-mp^b-p^b\sum_{\ell=1}^m j_\ell\right)(q^n+1)-iqM-kq^3 \right), j_1(q^n+1)+iM+k\right.\right.\] \[\left. ,\ldots,j_m(q^n+1)+iM+k  \Bigg); 1\leq k\leq M, 0\leq i\leq q, j_\ell\geq 0 \mbox{ and }\right.\]  \[\left. \left( q^2-mp^b-p^b\sum_{\ell=1}^m j_\ell\right)(q^n+1)-iqM-kq^3 >0  \right\}\;.\]
Then, $\Gamma(P_\infty,P_1,\ldots,P_m)=\Gamma_{m+1}$.
\end{theorem}
\begin{proof}
By Proposition \ref{Discrepancy}, we have that the divisor $A'$ is a discrepancy with respect to $P$ and $Q$ for any two distinct rational points $P,Q\in\chav{P_\infty,P_1,\ldots,P_m}$. By equivalence \eqref{effectivedivisor}, we can conclude that the divisor $A$ is also a discrepancy with respect to $P$ and $Q$ for any two distinct rational points $P,Q\in\chav{P_\infty,P_1,\ldots,P_m}$. Therefore, by Lemma \ref{LemaCT}, we have that $\Gamma_{m+1}\subseteq \Gamma(P_\infty,P_1,\ldots,P_m)$.

Next, we show that $\Gamma(P_\infty,P_1,\ldots,P_m)\subseteq \Gamma_{m+1}$. Let $\bold{n}=(n_0,n_1,\ldots,n_m)\in \Gamma(P_\infty,P_1\ldots,P_m)$. By Definition \ref{defnabla} and Proposition \ref{lubH}, we have that $\mathbf{n}$ is minimal in $\nabla_r(\mathbf{n})$ for all $0\leq r\leq m$. From Lemma \ref{mPoints}, $\mathbf{n}=(n_0,n_1,\ldots,n_m)\in G(P_\infty)\times G(P_1)\times\cdots\times G(P_m)$. Note that, as $H(P_\ell)=\langle q^n+1-iM-j: 0\leq i\leq p^b, 0\leq j\leq q^{n-3}(p^b-i)+i\dfrac{q^{n-3}-1}{q+1} \rangle$, for all $1\leq \ell\leq m$. Then, by the form of the elements in $G(P_\ell)$, we have that $n_\ell=j_\ell(q^n+1)+i_\ell M+k_\ell$, for some $j_\ell\geq 0, 0\leq i_\ell\leq q$ and $1\leq k_\ell\leq M$. Let
\[f=\dfrac{y^{q-i}z^{M-k}}{(x-\alpha_1)^{j_1+1}\cdots (x-\alpha_m)^{j_m+1}}\;,\] then 
\begin{eqnarray}
(f)_\infty & = &\dfrac{1}{p^b}\left( \left( q^2-mp^b-p^b\sum_{\ell=1}^m j_\ell\right)(q^n+1)-iqM-kq^3 \right)P_\infty \nonumber \\ & & +( j_1(q^n+1)+iM+k)P_1+\cdots+(j_m(q^n+1)+iM+k)P_m\;.\nonumber
\end{eqnarray}

We conclude that $f\in H(P_\infty,P_1,\ldots,P_m)$ and, as $(f)_\infty\sim A'$, then $(f)_\infty$ is a discrepancy with respect to $P$ and $Q$ for any rational points $P,Q\in \{P_\infty,P_1,\ldots,P_m\}$. So, by Lemma \ref{LemaCT}, we have that $\mathbf{f}=(\dfrac{1}{p^b}\left( \left( q^2-mp^b-p^b\sum_{\ell=1}^m j_\ell\right)(q^n+1)-iqM-kq^3 \right),$ $j_1(q^n+1)+iM+k,\ldots, j_m(q^n+1)+iM+k )\in \Gamma(P_\infty,P_1,\ldots,P_m)$.

Thus, $\mathbf{f}\in \nabla_r(\mathbf{n})$, for some $0\leq r\leq m$, and by Proposition \ref{minimal}, it follows that $\mathbf{f}$ is minimal in $\nabla_r(\mathbf{n})$ for all $r$, $0\leq r\leq m$. Furthermore, by minimality of $\mathbf{f}$ and $\mathbf{n}$, we have that $\mathbf{f}=\mathbf{n}$ and so $\Gamma(P_\infty,P_1,\ldots,P_m)\subseteq \Gamma_{m+1}$.
\end{proof}

\begin{example}
Using the values from Example \ref{ex1curve1}, we have the following divisors:
\begin{eqnarray}
(x-\alpha_\ell)&=&65P_\ell-65P_\infty\nonumber\\
(y)&=&13P_1+13P_2-26P_\infty\nonumber\\
(z)&=& P_1+P_2+Q_1+\cdots +Q_{30}-32P_\infty\nonumber
\end{eqnarray}
For this curve, taking $m=1$, by Theorem \ref{TeoGamma}, we have that $\Gamma(P_\infty,P_1)=\{(455-26i-65j-32k, 65j+13i+k):0\leq i\leq 4, j\geq 0, \mbox{ and }1\leq k\leq 13 \}$.

Taking $m=2$, we have that 
\begin{eqnarray}
\Gamma(P_\infty,P_1,P_2)&=&\{(390-26i-32k-65j_1-65j_2,65j_1+13i+k,65j_2+13i+k):\nonumber\\ & & 0\leq i\leq 4, j_1,j_2\geq 0, \mbox{ and }1\leq k\leq 13\}\;.\nonumber
\end{eqnarray}
\end{example}

\begin{example}
Using the values from Example \ref{ex2curve1}, we have the following divisors:
\begin{eqnarray}
(x-\alpha_\ell)&=& 1025P_\ell -1025P_\infty    \nonumber\\
(y)& = & 205P_1+205P_2 -410P_\infty  \nonumber \\
(z)& = & P_1+P_2+Q_1+\cdots +Q_{30} - 32P_\infty   \nonumber 
\end{eqnarray}
Taking $m=1$, by Theorem \ref{TeoGamma}, we have that $\Gamma(P_\infty,P_1)=\{(7175-410i-1025j-32k,205i+1025j+k):0\leq i\leq 4,j\geq 0 \mbox{ and } 1\leq k\leq 205\}$.

Taking $m=2$, we have that 
\begin{eqnarray}
\Gamma(P_\infty,P_1,P_2)& = &  \{(6150-410i-1025(j_1+j_2)-32k, 205i+1025j_1+k,205i+1025j_2+k): \nonumber \\
& & 0\leq i\leq 4,j_1,j_2\geq 0 \mbox{ and } 1\leq k\leq 205\}\; . \nonumber
\end{eqnarray}
\end{example}

\begin{example}
Using the values from Example \ref{ex1curve2}, we have the following divisors:
\begin{eqnarray}
(x-\alpha_\ell)&=&  28P_\ell -28P_\infty   \nonumber\\
(y)& = & 7P_1+7P_2+7P_3-21P_\infty   \nonumber \\
(z)& = & P_1+P_2+P_3+Q_1+\cdots +Q_{24}-27P_\infty   \nonumber 
\end{eqnarray}

For $m=1$, we have that 
\begin{eqnarray}
\Gamma(P_\infty,P_1)&=&  \{(224-21i-28j-27k,7i+28j+k): 0\leq i\leq 3, j\geq 0 \mbox{ and }1\leq k \leq 7\}\;. \nonumber 
\end{eqnarray}

For $m=2$, we have that 
\begin{eqnarray}
\Gamma(P_\infty,P_1,P_2)&=& \{(196-21i-28(j_1+j_2)-27k,7i+28j_1+k,7i+28j_2+k):   \nonumber \\
& & 0\leq i\leq 3, j_1,j_2\geq 0 \mbox{ and }1\leq k \leq 7\}\;.\nonumber
\end{eqnarray}

For $m=3$, we have that 
\begin{eqnarray}
\Gamma(P_\infty,P_1,P_2,P_3)&=& \{(168-21i-28(j_1+j_2+j_3)-27k,7i+28j_1+k,7i+28j_2+k,   \nonumber \\
& & 7i+28j_3+k): 0\leq i\leq 3, j_1,j_2,j_3\geq 0 \mbox{ and }1\leq k \leq 7\}\;.\nonumber
\end{eqnarray}

\end{example}

\begin{example}
Using the values from Example \ref{ex2curve2}, we have the following divisors:
\begin{eqnarray}
(x-\alpha_\ell)&=&  33P_\ell -33P_\infty   \nonumber\\
(y)& = & 11P_1+11P_2-22P_\infty   \nonumber \\
(z)& = & P_1+P_2+Q_1+\cdots +Q_6 -8P_\infty   \nonumber 
\end{eqnarray}

For $m=1$, we have that 
\begin{eqnarray}
\Gamma(P_\infty,P_1)&=&\{(99-22i-33j-8k,11i+33j+k):0\leq i\leq 2,j\geq 0 \mbox{ and } 1\leq k\leq 11\}\;.    \nonumber
\end{eqnarray}
\end{example}

For $m=2$, we have that 
\begin{eqnarray}
\Gamma(P_\infty,P_1,P_2)&=& \{(66-22i-33(j_1+j_2)-8k,11i+33j_1+k,11i+33j_2+k):    \nonumber \\
& & 0\leq i\leq 2,j_1,j_2\geq 0 \mbox{ and } 1\leq k\leq 11\}\;. \nonumber
\end{eqnarray}

\section{Pure Gaps and AG Codes}
In \cite{HK}, Homma and Kim introduced the concept of {\it pure gap}. An element $(n_1,\ldots,n_s)\in\mathbb{N}_0^s$ is a pure gap at $(P_1,\ldots,P_s)$ if 
\[ \ell\left( \sum_{i=1}^sn_iP_i-P_j \right)=\ell\left(
  \sum_{i=1}^sn_iP_i \right) \mbox{ for some }j\in\{1,\ldots,s\}\;.\]
Carvalho and Torres \cite[Lemma 2.5]{Carvalho-Torres} showed that
$(n_1,\ldots,n_s)$ is a pure gap at $(P_1,\ldots,P_s)$ if and only if
$\ell(\sum_{i=1}^sn_iP_i)=\ell(\sum_{i=1}^s(n_i-1)P_i)$. The authors
used this concept to obtain codes whose minimum distances have bounds
better than the Goppa bound. 

\begin{theorem}\cite[Theorem 3.3]{Carvalho-Torres}\label{codespuregaps}
Let $Q_1,\ldots,Q_n$, $P_1,\ldots,P_m$ be distinct $\mathbb{F}_q$-rational points of $\mathcal{X}$ and assume that $m\leq q$. Let $(\alpha_1,\ldots,\alpha_m),(\beta_1,\ldots,\beta_m)\in \mathbb{N}_0^m$ and set $D=Q_1+\cdots + Q_n$ and $G=\sum_{i=1}^m(\alpha_i+\beta_i-1)P_i$. Let $d_\Omega$ be the minimum distance of the code $C_{\Omega}(D,G)$. If $(\alpha_1,\ldots,\alpha_m),(\beta_1,\ldots,\beta_m)$ are pure gaps at $P_1,\ldots,P_m$, then $d_\Omega\geq \deg(G)-(2g-2)+m$, where $g$ is the genus of $\mathcal{X}$.
\end{theorem}

Using the same notation as in Section \ref{Weierstrass-S}, we
calculate the pure gaps at several points. The following results are stated in the same form as in \cite{CT}.

\begin{proposition}\cite[Proposition 4.2]{CT}\label{puregaps1}
Let $A=\sum_{\ell=0}^m a_\ell P_\ell$, where $(a_0,\ldots,a_m)\in \Gamma(P_0,\ldots,P_m)$. Let $\ell\in\{0,1,\ldots,m\}$, if $\mathcal{L}(A-P_\ell)=\mathcal{L}(A-2P_\ell)$, then $(a_0,a_1,\ldots,a_{\ell-1},a_\ell-1,a_{\ell+1}, \ldots,a_m)$ is a pure gap of $H(P_0,P_1,\ldots,P_m)$.
\end{proposition}


\begin{corollary}\cite[Corollary 4.3]{CT}\label{puregaps2}
If $2\leq k\leq M$, then $((q^2/p^b-m)(q^n+1)-kq^3,k,\ldots,k,k-1)$ is a pure gap of the Weierstrass semigroup $H(P_\infty,P_1,\ldots,P_m)$ on the $\mathcal{X}_{a,b,n,1}$ curve.

If $2\leq k\leq M$, then $((q^2-m)(q^n+1)-kq^3,k,\ldots,k,k-1)$ is a pure gap of the Weierstrass semigroup $H(P_\infty,P_1,\ldots,P_m)$ on the $\mathcal{Y}_{n,1}$ curve.
\end{corollary}

\begin{proposition}\cite[Proposition 4.4]{CT} \label{puregaps3}
Let $\alpha<2g-1$ and $(\alpha,1,\ldots,1)\in G(P_\infty,P_1,\ldots,P_m)$. If
\begin{enumerate}
	\item $\exists\, \lambda,\beta,\gamma\in\mathbb{N}_0$, with $\lambda\geq m$, such that $\lambda(q^n+1)+\beta qM+\gamma q^3=2g-1-\alpha$, or 
	\item $2g-2-\alpha\geq (m-1)(q^n+1)$ and $\exists\, \beta,\gamma\in\mathbb{N}_0$ such that $\beta qM+\gamma q^3=2g-1-\alpha$, 
\end{enumerate}
then $(\alpha,1,\ldots,1)$ is a pure gap.
\end{proposition}

In \cite{DAML}, the authors calculate the pure gaps in Kummer extensions defined by $y^m=f(x)$. The places $Q_1,\ldots,Q_r$ are all the zeros and poles of $f(x)$. They showed the following theorem, where $\lambda_i:=\upsilon_{Q_i}(f(x))$ denotes the multiplicity of the place $Q_i$. 

\begin{theorem}\cite[Theorem 3.3]{DAML}\label{TheoremPureGaps}
Let $P_1,\ldots,P_s\in \mathbb{P}_F$ be pairwise distinct totally ramified places in the Kummer extension $F/K(x)$. Then $(n_1,\ldots,n_s)\in\mathbb{N}_0^s$ is a pure gap at $(P_1,\ldots,P_s)$ if and only if for every $t\in \{0,\ldots,m-1\}$ exactly one of the two following conditions is satisfied:
\begin{enumerate}
	\item $\displaystyle\sum_{i=1}^s\left\lfloor
            \dfrac{n_i+t\lambda_i}{m} \right\rfloor+
          \displaystyle\sum_{i=s+1}^r\left\lfloor
            \dfrac{t\lambda_i}{m}  \right\rfloor < 0$ 
	
	\item $\left\lfloor \dfrac{n_i+t\lambda_i}{m}  \right\rfloor= \left\lfloor \dfrac{n_i-1+t\lambda_i}{m}  \right\rfloor$ for all $i\in \{1,\ldots,s\}$.
\end{enumerate} 
\end{theorem}

\begin{proposition}\cite[Proposition 3.9]{DAML}\label{prop-puregaps-Yn,s}
On the $\mathcal{Y}_{n,1}$ curve, let $P_1$ and $P_2$ be two totally ramified rational points that are different from $P_\infty$. Let $\alpha\in\{0,\ldots,q^2-3\}$ and $\beta\in\{0,1\}$. For $n\geq 5$, if
\[
\begin{array}{l}
n_1:= (\beta+1)q^{n-3}(q^2-q+1)+\alpha(q^n+1) \mbox{ and }\\
n_2:= (q^2-3)(q^n+1)+3q^{n-3}(q^2-q+1)-(\beta+1)q^{n-3}(q^2-q+1)-\alpha(q^n+1)\;,
\end{array}\]
 then the pair $(n_1,n_2)$ is a pure gap at $(P_1,P_2)$.
\end{proposition}

\begin{proposition}\cite[Proposition 3.10]{DAML}
On the $\mathcal{Y}_{n,1}$ curve, let $P_\infty$ be the unique rational point at infinity and $P_1$ be a totally ramified rational point different from $P_\infty$. For $\alpha\in\{0,\ldots,q^2-2\}$, the pair 
\[(n_1,n_2)=(1+\alpha(q^n+1),1+(q^2-2)(q^n+1)+q^n-2q^3+1-(1+\alpha(q^n+1)))\] is a pure gap at $(P_\infty,P_1)$.
\end{proposition}

\begin{proposition}
On the $\mathcal{X}_{a,b,n,1}$ curve, let $P_1$ and $P_2$ be two totally ramified rational points that are different from $P_\infty$. Let $\alpha\in\{0,\ldots,\frac{q^2}{p^b}-3\}$ and $\beta\in\{0,1\}$. For $n\geq 5$,~if  
\[
\begin{array}{l}
n_1:= (\beta+1)q^{n-3}(q^2-q+1)+\alpha(q^n+1) \mbox{ and }\\
n_2:= (\frac{q^2}{p^b}-3)(q^n+1)+3q^{n-3}(q^2-q+1)-(\beta+1)q^{n-3}(q^2-q+1)-\alpha(q^n+1)\;,
\end{array}\]  then the pair $(n_1,n_2)$ is a pure gap at $(P_1,P_2)$.

\end{proposition}

\begin{proof}
The rational points $P_1$ and $P_2$ are zeros of $t(x)$, and so $\lambda_1=\lambda_2=1$ in Theorem \ref{TheoremPureGaps}. Let $t\in\{0,\ldots,q^n\}$. We have that
$
\left\lfloor  \dfrac{n_1+t}{q^n+1}   \right\rfloor \neq \left\lfloor   \dfrac{n_1+t-1}{q^n+1}  \right\rfloor
$
 if and only if $n_1+t\equiv 0\pmod{q^n+1}$, which is equivalent to 
\[
\left\{ \begin{array}{ll}
t=q^n+1-(q^{n-1}-q^{n-2}+q^{n-3}) & \mbox{if } \beta=0 \\
t=q^n+1-(2q^{n-1}-2q^{n-2}+2q^{n-3}) & \mbox{if } \beta=1
\end{array} \right.
\]

Analogously,
$
\left\lfloor  \dfrac{n_2+t}{q^n+1}   \right\rfloor \neq \left\lfloor   \dfrac{n_2+t-1}{q^n+1}  \right\rfloor
$
 if and only if $n_2+t\equiv 0\pmod{q^n+1}$, which is equivalent to 
\[
\left\{ \begin{array}{ll}
t=q^n+1-(2q^{n-1}-2q^{n-2}+2q^{n-3}) & \mbox{if } \beta=0 \\
t=q^n+1-(q^{n-1}-q^{n-2}+q^{n-3}) & \mbox{if } \beta=1
\end{array} \right.
\]

Then we verified the first condition in Theorem \ref{TheoremPureGaps} for these values of $t$. Indeed, we have that 

\[
\left\lfloor  \dfrac{n_1+t}{q^n+1} \right\rfloor + \left\lfloor \dfrac{n_2+t}{q^n+1}  \right\rfloor +\dfrac{q}{p^b}(q-1)\left\lfloor \dfrac{t(q+1)}{q^n+1} \right\rfloor+\left\lfloor \dfrac{-tq^3/p^b}{q^n+1}  \right\rfloor = 
\]

\[\left\{
\begin{array}{ll}
\frac{q^2}{p^b}-1+\frac{q}{p^b}(q-1)(q-1)& \\ -\frac{q^3}{p^b}+\frac{1}{p^b}(q^2-q)=-1 & \mbox{if } t=q^n+1-(q^{n-1}-q^{n-2}+q^{n-3})\\
\frac{q^2}{p^b}-2+\frac{q}{p^b}(q-1)(q-2)& \\ -\frac{q^3}{p^b}+\frac{1}{p^b}(2q^2-2q)+1=-1 & \mbox{if } t=q^n+1-(2q^{n-1}-2q^{n-2}+2q^{n-3})\;.
\end{array}\right.
\]
\end{proof}


\begin{remark}
The parameters of the AG codes over the curves $\mathcal{X}_{a,b,n,1}$ and $\mathcal{Y}_{n,1}$ cannot be compared with the parameters of the codes in MinT's tables \cite{MinT} since the size of the corresponding alphabet is too large. However, the relative parameters of these codes can be compared with the relative parameters of the AG codes constructed from the GGS curves or induced by them. Given a $[n,k,d]_{\mathbb{F}_q}$ linear code, we have that the relative parameters are $k/n$ the rate and $d/n$ the relative minimum distance, and by the Singleton bound we have that $\frac{k}{n}+\frac{d}{n}\leq 1 +\frac{1}{n}$. 
\end{remark}

\begin{example}
Consider the curve $\mathcal{Y}_{5,1}$ in Example \ref{ex2curve2} over $\mathbb{F}_{2^{10}}$. By Proposition \ref{prop-puregaps-Yn,s}, taking $\alpha=1$ we have that $(34,50)$ is a pure gap at $P_\infty,P_1$. By Theorem \ref{codespuregaps}, we have that the two-point code $C_\Omega(D, 67P_\infty+99P_1)$ has minimum distance $d_\Omega\geq 78$, hence yielding a $[3967,3846,\geq 78]_{2^{10}}$ code. This code has better relative parameters than the corresponding one-point AG code $[3968,3846,\geq 77]_{2^{10}}$ given in \cite[Table 1]{BMZ1}. 
\end{example}

\begin{example}
Consider the curve $\mathcal{X}_{2,1,3,1}$ over $\mathbb{F}_{4^6}$ given in Example \ref{ex1curve1}. By Proposition \ref{puregaps3}, it follows that $(230,1)$ is a pure gap at $P_\infty,P_1$. By Theorem \ref{codespuregaps}, we have that the two-point code $C_\Omega(D, 459P_\infty+P_1)$ 
has minimum distance $d_\Omega\geq 40$,
hence yielding a $[n=31231,k=30982,d\geq 40]$ code. 
The bound on the minimum distance is better than the one corresponding to the one-point AG code given in \cite[Corollary 5.5 (2)]{TTT}, defined over the same curve, whose Feng-Rao bound for the minimum distance is $\delta_{FR}(249)=39$, hence yielding a
$[n=31232,k=30982,d\geq 39]$ code.
\end{example}

\begin{example}
Consider the curve $\mathcal{X}_{1,1,3,1}$ of genus $g=3$ over $\mathbb{F}_{2^6}$. We have that $H(P_\infty)=\langle 3,4 \rangle$ and as the pole divisor $(z/y^2)_\infty=5P_0$, $(z^2/y^3)_\infty=7P_0$ then $H(P_0)=\langle 3,5,7 \rangle$. For this curve, by Theorem \ref{TeoGamma}, we have that $\Gamma(P_\infty,P_0)=\{(5,1),(1,2),(2,4)\}$.
Take the divisor $G=4P_\infty+P_0$. 
Using the {\sc Magma} software, 
one can see that the two-point AG code $C_\Omega(D,4P_\infty+P_0)$ has parameters $[111,108,3]$. This code is quasi perfect.
\end{example}

\bibliographystyle{spmpsci}



\appendix

\section{Genus of $S$}
\label{a:S}
In this section we will 
determine the genus of the semigroup $S$
generated by 
  $\{q^3+1-iN-j:0\leq i+j\leq p^b\},$
  with $q=(p^b)^r$, for a prime number $p$, some positives integers $b,r$, with $r\geq 2$, and $N=\frac{q^3+1}{q+1}$. 

For this purpose 
we will first consider 
the numerical semigroup $S_{P,N,K}$ generated by 
$\{KN+aN-j:0\leq j\leq a\leq P\}$,
where 
$P$, $N$, $K$ are positive integers with
$P\mid N-1$,
$P\mid K-1$,
$K<N$.
Notice that $S=S_{P,N,K}$ with $K=q+1-p^b$ and $P=p^b$, with $a$
playing the role of $p^b-i$, and $j,N$ playing their own role. The required conditions hold, indeed, 
  $P\mid N-1$ since $N-1=\frac{q^3+1}{q+1}-1=(q^2-q+1)-1=q^2-q$ and, similarly,
  $P\mid K-1$.

\subsection{Characterization}

With the same notation as before, for an integer $M\geq 0$, let $S^M_{P,N,K}=\{MKN+aN-j:0\leq j\leq
a\leq MP\}$ and let $\widetilde S^M_{P,N,K}=\{MKN+aN-j: \max\{0,(M-1)P-K+1\}\leq a \leq MP,  0\leq j\leq \min\{a,N-1\}\}$.
It is obvious that 
$S_{P,N,K}=\cup_{M\geq 0} S^M_{P,N,K}$ and that $\widetilde
S^M_{P,N,K}\subseteq S^M_{P,N,K}$.
Now, 
any 
element in $S^M_{P,N,K}$ is in at least one set $\widetilde S^{M'}_{P,N,K}$ for 
some $M'\leq M$. This can be proved by induction on $M$. For $M=0$ and for $M=1$ it
is straightforward.
For $M>1$, suppose that an element of $S^M_{P,N,K}$ is $\ell=MKN+aN-j$ for 
some particular $0\leq j\leq a\leq MP$. If $a\geq (M-1)P-K+1$ and 
$j\leq N-1$ then $\ell\in\widetilde S^M_{P,N,K}$. Otherwise, if 
$0\leq a<(M-1)P-K+1$,
then $\ell=(M-1)KN+(K+a)N-j=(M-1)KN+a'-j$ with $a'=K+a\leq (M-1)P$ and 
$0\leq j\leq a\leq a'$, so $\ell\in S^{M-1}_{P,N,K}$ and the result
follows by induction.
If $a\geq (M-1)P-K+1$ but  
$j>N-1$, then suppose that $Q$ is the quotient of the division of 
$j$ by $N$. Then 
$\ell=MKN+(a-Q)N-(j-QN)=MKN+a'-j'$ with $a'=a-Q \geq a-j \geq 0$, and
$a'\leq a\leq MP$. Furthermore,
$j'=j-QN$, which is the remainder of the division of 
$j$ by $N$, and which is between $0$ and $N-1$. So, $\ell\in \widetilde S^M_{P,N,K}$.
Consequently, we also have $S_{P,N,K}=\cup_{M\geq 0} \widetilde S^M_{P,N,K}$.

Now, $\widetilde S^M_{P,N,K}\subseteq [(M-1)(K+P)N+1,M(K+P)N]$. Indeed, the
minimum of $\widetilde S^M_{P,N,K}$ is at least $(M-1)(KN+PN)+1$ since the elements in $\widetilde S^M_{P,N,K}$ satisfy
$MKN+aN-j\geq MKN+((M-1)P-K+1)N-N+1=MKN+MPN-PN-KN+N-N+1=(M-1)(K+P)N+1$.
On the other hand, the maximum of $\widetilde S^M_{P,N,K}$ is at most $M(K+P)N$
since the elements in $\widetilde S^M_{P,N,K}$ satisfy
$MKN+aN-j\leq MKN+MPN-0=M(K+P)N$.

In particular, the sets $\widetilde{S}^M_{P,N,K}$ are disjoint and, so,
$S_{P,N,K}=\sqcup_{M\geq 0}\widetilde S^M_{P,N,K}.$

Let $M_0=\frac{K-1}{P}+1$, $M_1=\frac{N-1}{P}$, $M_2=\frac{K+N-2}{P}$. 

\subsection{Conductor of $S_{P,N,K}$}

Now we are ready to determine the 
Frobenius number of $S_{P,N,K}$, that is, its largest gap.
The conductor of $S_{P,N,K}$ is then the non-gap of $S_{P,N,K}$ right after its
Frobenius number.

Observe that $\widetilde{S}_M=\cup_{a\geq 
  \max\{0,(M-1)P-K+1\}}I_a$ with 
$I_a=\{MKN+aN-j: 0\leq j\leq \min\{a,N-1\}\}$. 
Since $j$ ranges from $0$ to $\min\{a,N-1\}$, there are gaps between the intervals $I_a$ and $I_{a-1}$ if and only 
if $\min\{a,N-1\}<N-1$, i.e., if and only if $a< N-1$. Since in
$\widetilde S^M_{P,N,K}$ 
$a$ ranges from $\max\{0,(M-1)P-K+1\}$ to $PN$, the inequality $a<N-1$
occurs in $\widetilde S^M_{P,N,K}$ if and only if $(M-1)P-K+1<N-1$, that is, if and only if $M\leq\frac{N+K-2}{P}$.

Let $M_F= \frac{N+K-2}{P}$.
The last gap of $S_{P,N,K}$ will then be the gap previous to $I_a$
with $a=N-2$ in $\widetilde S^{M_F}_{P,N,K}$, that is, the Frobenius
number will be $M_FKN+aN-a-1$ for $a=N-2$, i.e. 
$\frac{K+N-2}{P}KN+(N-2)N-N+1=\frac{K+N-2}{P}KN+N^2-3N+1$.

Simplifying by means of {\sc Sage} we obtain that the conductor of
$S_{P,N,K}$ is $$c=\frac{N^2K + NK^2 + N^2P - 2NK - 3NP + 2P}{P}.$$

\subsection{Genus}

Let $M_0=\frac{K-1}{P}+1$, $M_1=\frac{N-1}{P}$, $M_2=\frac{K+N-2}{P}=M_0+M_1-1$. 
  \begin{itemize}
  \item   If $0\leq M\leq M_0$ then   $S_{P,N,K}^M=\widetilde S_{P,N,K}^M$ and $\#\widetilde S_{P,N,K}^M=\sum_{a=0}^{MP}(a+1)=\sum_{b=1}^{MP+1} b=\frac{(MP+1)(MP+2)}{2}$.
  \item   If $M_0<M\leq M_1$ then $\#\widetilde S_{P,N,K}^M=\sum_{a=(M-1)P-K+1}^{MP}(a+1)=\sum_{(M-1)P-K+2}^{MP+1}b=\frac{(MP+1)(MP+2)}{2}-\frac{(MP-P-K+1)(MP-P-K+2)}{2}$. This is because if $M\leq M_1$, then $MP\leq N-1$, so $j$ will always range from $0$ to $a$.
  \item If $M>M_1$ then $\#([(M-1)(K+P)N+1,M(K+P)N]\setminus
    \widetilde{S}_{P,N,K}^M)=$\break$\sum_{a=(M-1)P-K+1}^{N-2}(N-a-1)=\sum_{b=1}^{N-(M-1)P+K-2}b=$\break$\frac{(N-MP+P+K-2)(N-MP+P+K-1)}{2}$. This
    is because if $M>M_1$ then $MP>N-1$ and so at some point $a=N-1$. In this case, $[(M-1)(K+P)N+1,M(K+P)N]
\setminus\widetilde S_{P,N,K}^M=\{MKN+aN-j: (M-1)P-K+1 \leq a \leq N-2,
a+1\leq j\leq N-1\}$.
\end{itemize}
Now, using the formulas 
$\sum_{M=1}^n M=\frac{n(n+1)}{2}$
and $\sum_{M=1}^n M^2=\frac{n(n+1)(2n+1)}{6}$
we get to the final count of the genus:
\begin{eqnarray*}
g&=&(M_1(P+K)N)\\
&&-\sum_{M=1}^{M_1}\frac{(MP+1)(MP+2)}{2}\\
&&+\sum_{M=M_0+1}^{M_1}\frac{(MP-P-K+1)(MP-P-K+2)}{2}\\
&&+\sum_{M=M_1+1}^{M_2}\frac{(N-MP+P+K-2)(N-MP+P+K-1)}{2}\\
&=&(M_1(P+K)N)\\
&&-\frac{1}{12}P^2M_1(M_1+1)(2M_1+1)-\frac{3}{4}PM_1(M_1+1)-M_1\\
&&+\frac{1}{12}P^2(M_1(M_1+1)(2M_1+1)-M_0(M_0+1)(2M_0+1))\\&&+\frac{1}{4}P(3-2P-2K)(M_1(M_1+1)-M_0(M_0+1))\\&&+\frac{1}{2}(P+K-1)(P+K-2)(M_1-M_0)\\
&&+\frac{1}{12}P^2(M_2(M_2+1)(2M_2+1)-M_1(M_1+1)(2M_1+1))\\&&+\frac{1}{4}P(-2(P+K+N)+3)(M_2(M_2+1)-M_1(M_1+1))\\&&+\frac{1}{2}(P+K+N-1)(P+K+N-2)(M_2-M_1)\\
&=&\frac{N^2K + NK^2 + N^2P + NKP - 3NK - 3NP + P + 1}{2P}
\end{eqnarray*}

\subsection{Back to the originary problem}
If we take $N=\frac{q^3+1}{q+1}$, $K=q+1-p^b$, $P=p^b$, then the genus
is $g=\frac{q^5 - q^3p^b - q^3 + q^2}{2p^b},$
while the conductor is
$c=\frac{q^5 - 2q^3p^b + q^2p^{2b} - q^2p^b - qp^{2b} + q^2 + qp^b + p^{2b} - p^b}{p^b}.$

\section{Genus of $S'$}
\label{a:Sprime}
Suppose we have
$q=p^a$, $b\mid a$, $b\neq a$,
$n$ odd, $n\geq 3$.
Let $M=\frac{q^n+1}{q+1}=q^{n-1}-q^{n-2}+q^{n-3}-q^{n-4}+\dots-q+1$.
We want to prove that the genus of the semigroup
$S'=\langle q^n+1-iM-j:0\leq iM+jq^2\leq q^{n-1}p^b \rangle$
is 
$\frac{q^{n+2}-p^bq^n-q^3+q^2}{2p^b}$.


\subsection{Definition of $S'$ revisited}

\begin{lemma}\label{l:bone}
  If $n\geq 3$,
$$S'=\left\langle k(q^{n-1}-q^{n-2})+\ell
:q+1-p^b\leq k\leq q+1 \mbox{ and }
q^{n-3}(q-p^b) \leq \ell \leq k\frac{q^{n-2}+1}{q+1}
\right\rangle$$
Equivalently, 
by setting
$A=q^{n-1}-q^{n-2}$,
$k_0=q+1-p^b$,
$k_1=q+1$,
$\ell_0=q^{n-3}(q-p^b)+1$,
$\ell_1=\frac{q^{n-2}+1}{q+1}$,
the semigroup $S'$ is 
$S'=\left\langle kA+\ell: 
k_0\leq k\leq k_1,
\ell_0 \leq \ell \leq k\ell_1\right\rangle.$
\end{lemma}

\begin{proof}
  We can rewrite $S'$ as
  $S'=\langle (q+1-i)M-j:0\leq iM+jq^2\leq q^{n-1}p^b \rangle.$ The integer $i$ is then bounded as $0\leq i\leq \lfloor\frac{q^{n-1}p^b}{M}\rfloor$. The quotient and the remainder of the division of $q^{n-1}p^b$ by $M$ are, respectively, $p^b$ and 
$p^b(q^{n-2}-q^{n-3}+\dots+q-1)$ (since this remainder is between $0$ and $M-1$).
Consequently, 
$0\leq i\leq p^b$.
Now, setting $k=q+1-i$, the bounds for $k$ are 
$q-p^b+1\leq k\leq q+1$.
Hence, since $iM=(q+1-k)M=(q+1)M-kM=q^n+1-kM$,
$$S'=\langle kM-j:q-p^b+1\leq k\leq q+1\mbox{ and }0\leq q^n+1-kM+jq^2\leq q^{n-1}p^b \rangle$$

Finally, we want to replace the bounds for $q^n+1-kM+jq^2$ by bounds on $j$. Reorganizing them, we obtain
$$kM-q^n-1\leq jq^2\leq kM-q^n+q^{n-1}p^b-1=kM-q^{n-1}(q-p^b)-1.$$
Since $k\leq q+1$, the lower bound is non-positive. So,
$0\leq jq^2$

As for the upper bound on $j$, $j\leq
\left\lfloor\frac{kM-q^{n-1}(q-p^b)-1}{q^2}\right\rfloor=$ \break $
\left\lfloor\frac{kq^{n-1}-kq^{n-2}+kq^{n-3}-\dots+kq^2-kq+k-q^{n-1}(q-p^b)-1}{q^2}\right\rfloor=
\left\lfloor k\frac{q^{n-2}+1}{q+1}+\frac{-kq+k-1}{q^2}-q^{n-3}(q-p^b)\right\rfloor
=\left\lfloor k\frac{q^{n-2}+1}{q+1}
-q^{n-3}(q-p^b)-1
+\frac{q^2-k(q-1)-1}{q^2}
\right\rfloor$.
By the bounds on $k$ we deduce that 
$0\leq \frac{q^2-k(q-1)-1}{q^2}\leq \frac{p^b(q-1)}{q^2}<1$. So,
$$S'=\langle kM-j:q+1-p^b\leq k\leq q+1 \mbox{ and }
0\leq j\leq k\frac{q^{n-2}+1}{q+1}-q^{n-3}(q-p^b)-1 
\rangle$$

Let now $\ell=k\frac{q^{n-2}+1}{q+1}-j$. Notice that 
$kM-j=k(M-\frac{q^{n-2}+1}{q+1})+\ell
=k(q^{n-1}-q^{n-2})+\ell$.
The bounds of $\ell$ are
$q^{n-3}(q-p^b)+1 \leq \ell \leq k\frac{q^{n-2}+1}{q+1}$.
\end{proof}

Let $G=\left\{ kA+\ell: 
k_0\leq k\leq k_1,
\ell_0 \leq \ell \leq k\ell_1
\right\}$ and let $mG=\{a_1+\dots+a_m:a_i\in G\}$.
Define $B_m=[(m-1)(q^n+1)+1,m(q^n+1)]\cap mG$.
\begin{lemma}\label{l:btwo}
The following statements hold.
\begin{enumerate}
\item $S'=\{0\}\cup\bigcup_{m\geq 1} mG$,
  \item     
$S'=\{0\}\cup \sqcup_{m\geq 1} B_m$.
  \end{enumerate}
\end{lemma}

\begin{proof}
  First of all, notice that 
$mG=\underbrace{[mk_0A+m\ell_0,mk_0(A+\ell_1)]}_{mk_0\ell_1 -m\ell_0
  +1 } $
$\cup$\break
        $\underbrace{[(mk_0+1)A+m\ell_0,(mk_0+1)(A+\ell_1)]}_{(mk_0+1)\ell_1
          -m\ell_0 +1 }$
$\cup
        \dots
\cup$
        $\underbrace{[mk_1A+m\ell_0,mk_1(A+\ell_1)]}_{mk_1\ell_1-m\ell_0+1}$
        \begin{enumerate}
        \item
          The first part is obvious and follows from the definitions.
          \item For the second part, it is obvious that the sets $B_m$ are disjoint and it is obvious the inclusion $\supseteq$.
  Let us prove for all $m$ the inclusion $S'\cap[(m-1)(q^n+1)+1,m(q^n+1)]\subseteq
  B_m$ by induction on $m$.

  First of all we need to see that
  $S'\cap[1,q^n+1]=B_1$.
  The smallest element of $2G$ is $2(k_0A+\ell_0)=
  2((q+1-p^b)(q^{n-1}-q^{n-2})+q^{n-3}(q-p^b)+1)=
  2(q^n+1-p^b(q^{n-1}-q^{n-2}+q^{n-3}))=
  q^n+1+(q^n+1-2p^b(q^{n-1}-q^{n-2}+q^{n-3}))
  >
  q^n+1+(q^n+1-2p^b\frac{q^n+1}{q+1})
  \geq q^n+1$ if $2p^b\leq q+1$, which is a consequence of the fact that $p^b<q$.
  
Now suppose $m>1$.
Since the maximum of $mG$ is $m(q^n+1)$,
  we have $mG\subseteq [0,m(q^n+1)]$.
  Now it will suffice to see that $mG\cap [0,(m-1)(q^n+1)]\subseteq (m-1)G$ and the result will follow by induction.

  Notice that $mG$ is the union of the sets of the form $S_{m,\tilde k}=[\tilde k A+m\ell_0,\tilde k(A+\ell_1)]$ for some $\tilde k$ satisfying $mk_0\leq \tilde k\leq m k_1$, while 
  $(m-1)G$ is the union of sets of the form $S_{(m-1),\tilde{\tilde{k}}}[\tilde{\tilde{k}} A+m\ell_0,\tilde{\tilde{k}}(A+\ell_1)]$ for some $\tilde{\tilde{k}}$ satisfying $(m-1)k_0\leq \tilde{\tilde{k}}\leq (m-1) k_1$.

  Suppose that $a\in mG\cap [0,(m-1)(q^n+1)]$.
  If $a\in S_{m\tilde k}$ with $mk_0\leq \tilde k\leq (m-1)k_1$, then, since $S_{m,\tilde k}\subseteq S_{m-1,\tilde k}$, we have $a\in S_{m-1,\tilde k}\subseteq (m-1)G$.
  On the other hand, if $a\in S_{m,\tilde k}\cap [0,(m-1)(q^n+1)]$ with $(m-1)k_1<\tilde k \leq mk_1$,
then $a\geq \tilde k A+ m\ell_0 >(m-1)k_1 A +(m-1)\ell_0$.
So, $a\in S_{m-1,(m-1)k_1}\subseteq (m-1)G$.
\end{enumerate}
  \end{proof}

\subsection{Number of gaps by intervals}

Let $C_m=[(m-1)(q^n+1)+1,m(q^n+1)]\setminus B_m$.
In this section we wonder what are the elements in $C_m$.
As before, we split the elements in $mG$ into (not necessarily disjoint) blocks
of the form
$S_{m,\tilde k}=[\tilde k A+m\ell_0,\tilde k(A+\ell_1)]$ for some $\tilde k$ satisfying $mk_0\leq \tilde k\leq m k_1$.


\begin{lemma}\label{l:bthree}
\begin{enumerate}
\item 
  Suppose that $mk_0< k\leq mk_1$. Then 
the gaps between $S_{m,(k-1)}$ and $S_{m,k}$
are contained in 
$[(m-1)(q^n+1)+1,m(q^n+1)]$ if and only if 
$ k\geq \max\left(mk_0+1,mq+m-q\right)$
\item 
$\max\left(mk_0+1,mq+m-q\right)=mk_0+1$ if and only if $m\leq M_1:=p^{a-b}$. 
\end{enumerate}
\end{lemma}

\begin{proof}
\begin{enumerate}
\item 
   Suppose that $mk_0<  k\leq mk_1$. Then 
the gaps between $S_{m,(k-1)}$ and $S_{m,k}$
are contained in 
$[(m-1)(q^n+1)+1,m(q^n+1)]$ if and only if 
$(k-1)\frac{q^n+1}{q+1}\geq (m-1)(q^n+1)$, that is, if and only if 
$k\geq (m-1)(q+1)+1=mq+m-q$. 
\item 
$\max\left(mk_0+1,mq+m-q\right)=mk_0+1$ if and only if 
$mq+m-q\leq m(q+1-p^b)+1$, that is, if and only if 
$-q\leq -mp^b+1$, i.e., $mp^b\leq q+1$. 
Now observe that the quotient of the Euclidean division of $q+1$ by 
$p^b$
is $p^{a-b}$ while the remainder is $1$. 
So, the statement follows. 
\end{enumerate}
\end{proof}

\begin{lemma} \label{l:bfour}
\begin{enumerate}
\item Suppose that $mk_0< k\leq mk_1$. Then
  there are gaps between
  $S_{m,(k-1)}$ and $S_{m,k}$ if and only if 
  $k\leq \min\left(
  mk_1,\frac{q^{n}-q-1+m(q+1)(q^{n-3}(q-p^b)+1)}{q^{n-2}+1}\right).$
\item If $n>3$,
$\min\left(mk_1,
\frac{q^{n}-q+m(q+1)(q^{n-3}(q-p^b)+1)}{q^{n-2}+1}\right)=mk_1$
if and only if $m\leq M_2:=(q-1)p^{a-b}$.
\item 
If $n=3$,
$\min\left(mk_1,
\frac{q^{n}-q+m(q+1)(q^{n-3}(q-p^b)+1)}{q^{n-2}+1}\right)=mk_1$
if and only if $m\leq \tilde{M_2}:=(q-1)p^{a-b}-1$.
\end{enumerate}
\end{lemma}

\begin{proof}
\begin{enumerate}
\item
Suppose that $mk_0< k\leq mk_1$. Then there are gaps
between $S_{m,( k-1)}$ and $S_{m, k}$ if and only if 
$( k-1)(A+\ell_1)\leq k A+m\ell_0-2$, equivalently,
$( k-1)\ell_1\leq m\ell_0-2+A$, equivalently,
$ k\leq\frac{m\ell_0-2+A}{\ell_1}+1=
(q+1)\frac{m(q^{n-3}(q-p^b)+1)-2+q^{n-1}-q^{n-2}}{q^{n-2}+1}+1=
\frac{m(q+1)(q^{n-3}(q-p^b)+1)-q-2+q^{n}-q^{n-2}+q^{n-2}+1}
{q^{n-2}+1}
=
\frac{
q^{n}
-q-1
+m((q+1)q^{n-3}(q-p^b)+1)
}{q^{n-2}+1}$.
\item
$\min\left(mk_1,
\frac{q^{n}-q-1+m(q+1)(q^{n-3}(q-p^b)+1)}{q^{n-2}+1}\right)=mk_1$
if and only if
$m(q+1)\leq \frac{q^{n}-q-1+m(q+1)(q^{n-3}(q-p^b)+1)}{q^{n-2}+1}$, that is,
if and only if
$m(q^{n-1}+q^{n-2}+q+1)\leq q^{n}-q-1+m(q+1)(q^{n-3}(q-p^b)+1)$, i.e.,
$m(q^{n-1}+q^{n-2}+q+1)\leq q^{n}-q-1+m(q^{n-2}(q-p^b)+q^{n-3}(q-p^b)+q+1)$, i.e.,
$0\leq q^{n}-q-1-mp^b(q+1)q^{n-3}$, i.e.,
$m\leq \lfloor\frac{q^n-q-1}{p^b(q+1)q^{n-3}}\rfloor$

Here we notice that the Euclidean division of $q^n-q-1$ by 
$q^{n-2}p^b+q^{n-3}p^b$ has quotient $(q-1)p^{a-b}$ and remainder
$q^n-q-1-(q-1)p^{a-b}(q^{n-2}p^b+q^{n-3}p^b)=
q^n-q-1
-q^{n}-q^{n-1}
+q^{n-1}+q^{n-2}
=
q^{n-2}
-q-1
$.

So, the statement follows.
\item It can be proved as the previous item.
\end{enumerate}
\end{proof}

\begin{lemma}\label{l:bfive}
If $(q-1)p^{a-b}+1\leq m\leq qp^{a-b}-1$ then
$\frac{q^{n}-q-1+m(q+1)(q^{n-3}(q-p^b)+1)}{q^{n-2}+1}$
is not an integer and
$\lfloor\frac{q^{n}-q-1+m(q+1)(q^{n-3}(q-p^b)+1)}{q^{n-2}+1}
\rfloor=q^2-q+m(q-p^b+1)$.
\end{lemma}

\begin{proof}
$\frac{q^{n}-q-1+m(q+1)(q^{n-3}(q-p^b)+1)}{q^{n-2}+1}=
\frac{q^{n}-q+1}{q^{n-2}+1}+
m\frac{(q+1)(q^{n-3}(q-p^b)+1)}{q^{n-2}+1}
=q^2-\frac{q^2+q+1}{q^{n-2}+1}
+m(q-p^b+1)-mp^b\frac{q^{n-3}-1}{q^{n-2}+1}=
q^2-q+m(q-p^b+1)+q-\frac{q^2+q+1}{q^{n-2}+1}-mp^b\frac{q^{n-3}-1}{q^{n-2}+1}.
$

Now, it is enough to see that 
$0\leq q-\frac{q^2+q+1}{q^{n-2}+1}-mp^b\frac{q^{n-3}-1}{q^{n-2}+1}<1$.

On one hand, 
$
q-\frac{q^2+q+1}{q^{n-2}+1} -
mp^b(
\frac{q^{n-3}-1}{q^{n-2}+1}) 
\geq 
q-
\frac{q^2+q+1}{q^{n-2}+1} 
-(qp^{a-b}-1)(\frac{(q^{n-3}-1)p^b}{q^{n-2}+1}) 
=
q-
\frac{q^2+q+1}{q^{n-2}+1} -
\frac{(q^{n-3}-1)q^2-p^{b}(q^{n-3}-1)}{q^{n-2}+1}
=
q- 
\frac{q^2+q+1}{q^{n-2}+1} -
\frac{q^{n-1}-q^2-q^{n-3}p^b+p^b}{q^{n-2}+1}
=
q+
\frac{-q^{n-1}+q^{n-3}p^b-p^b+q+1}{q^{n-2}+1}
=
\frac{q+q^{n-3}p^b-p^b+q+1}{q^{n-2}+1}
=
\frac{2q+p^b(q^{n-3}-1)+1}{q^{n-2}+1}
>0$.

On the other hand,
$
q-
\frac{q^2+q+1}{q^{n-2}+1} -
mp^b\frac{q^{n-3}-1}{q^{n-2}+1}
\leq 
\frac{q^{n-1}+q-q^2-q-1-(qp^{a-b}-p^{a-b}+1)(q^{n-3}p^b-p^b) }{q^{n-2}+1}=
\frac{q^{n-1}-q^2-1 
-qp^{a-b}(q^{n-3}p^b-p^b) 
+p^{a-b}(q^{n-3}p^b-p^b) 
-(q^{n-3}p^b-p^b) 
}{q^{n-2}+1}
=
\frac{q^{n-1}-q^2-1 
-q^{n-1}+q^2
+q^{n-2}-q
-q^{n-3}p^b+p^b 
}{q^{n-2}+1}
=
\frac{
-1 
+q^{n-2}-q
-q^{n-3}p^b+p^b 
}{q^{n-2}+1}
=
\frac{
  q^{n-2}+1-(2+q+(q^{n-3}-1)p^b)
}{q^{n-2}+1}
<1.
$

\end{proof}

\begin{lemma}\label{l:bsix}
     Suppose that $mk_0< k\leq mk_1$. Then
the number of gaps between $S_{m,(k-1)}$ and $S_{m,k}$
is
$m(q^{n-3}(q-p^b)+1)-{ k}\frac{q^{n-2}+1}{q+1}+\frac{q^n+1}{q+1}-1$.
\end{lemma}

\begin{proof}
The number of gaps between $S_{m,(k-1)}$ and $S_{m, k}$
is
$kA+m\ell_0-(k-1)(A+\ell_1)-1=
m\ell_0-k\ell_1 +A+\ell_1-1$, which yields the formula in the statement. \end{proof}

\begin{lemma}\label{l:bseven}
\begin{enumerate}
\item
There are gaps that are at least $(m-1)(q^n+1)+1$
and which are smaller than the elements in $S_{m,mk_0}$ if and only if
$m\leq p^{a-b}=M_1$. 
\item
If $m\leq p^{a-b}=M_1$,
then the number of gaps between $(m-1)(q^n+1)+1$
  and $S_{m,mk_0}$
  is 
  $q^n-mp^b(q^{n-1}-q^{n-2}+q^{n-3})$.
\end{enumerate}
  \end{lemma}

\begin{proof}
\begin{enumerate}
\item
There are gaps that are at least $(m-1)(q^n+1)+1$
and which are smaller than the elements in $S_{m,mk_0}$ if and only if 
   $mk_0A+m\ell_0\geq (m-1)(q^n+1)+2$. This is equivalent to
    $m(q+1-p^b)(q^{n-1}-q^{n-2})+m (q^{n-3}(q-p^b)+1)\geq (m-1)(q^n+1)+2$, that is,
if and only if  
$
mq^n-mq^{n-1}+mq^{n-1}-mq^{n-2}-mq^{n-1}p^b+mq^{n-2}p^b 
+
mq^{n-2}-mq^{n-3}p^b+m 
\geq 
mq^n 
+m 
-q^n-1 
+2$,
which is equivalent to 
$q^n-1\geq mp^b(q^{n-1}-q^{n-2}+q^{n-3})$ i.e., 
$m\leq \lfloor\frac{q^n-1}{p^b(q^{n-1}-q^{n-2}+q^{n-3})}\rfloor$. 

Here we remark that the Euclidean division of $q^n-1$ by 
 $p^b(q^{n-1}-q^{n-2}+q^{n-3})$
has quotient $p^{a-b}$ and remainder 
$q^n-1-(q^{n}-q^{n-1}+q^{n-2}) 
=q^{n-1}-q^{n-2}-1$.
So, the result follows.
\item
  It follows from the formula
$mk_0A+m\ell_0-(m-1)(q^n+1)-1$
and a similar simplification as before.
\end{enumerate}
\end{proof}

\begin{lemma}\label{l:beight}
  Let $M_3=qp^{a-b}-1$. The set $C_m$ is not empty if and only if $m\leq M_3$.
  \end{lemma}
\begin{proof}
  It is clear that for $m<M_2$, $C_m\neq \emptyset$.
  For $m\geq M_2$, $C_m\neq \emptyset$ if and only if
  $$\frac{q^n-q-1+m(q+1)(q^{n-3}(q-p^b)+1)}{q^{n-2}+1}\geq m(q+1)-q.$$
  This is equivalent to 
  $q^n-q-1+m(q+1)(q^{n-3}(q-p^b)+1-(q^{n-2}+1))\geq -q(q^{n-2}+1)$, which in turn is equivalent to $q^n-1-m(q+1)(p^bq^{n-3})\geq -q^{n-1}$, i.e.,  $m\leq \lfloor\frac{(q+1)q^{n-1}-1}{(q+1)p^bq^{n-3}}\rfloor=qp^{a-b}-1$. 
\end{proof}

\begin{corollary}\label{c:bnine}
$S'=\sqcup_{m=1}^{M_3}B_m$
\end{corollary}


\begin{theorem}
The genus of $S'$ is $\frac{q^{n+2}-p^bq^n-q^3+q^2}{2p^b}$.
\end{theorem}

\begin{proof} By Lemma~\ref{l:bone}, Lemma~\ref{l:btwo}, Lemma~\ref{l:bthree}, Lemma~\ref{l:bfour}, Lemma~\ref{l:bfive}, Lemma~\ref{l:bsix}, Lemma~\ref{l:bseven}, Lemma~\ref{l:beight}, and Corollary~\ref{c:bnine}, it easily follows that the genus of $S'$ is
\begin{eqnarray*}
&&\sum_{m=1}^{M_1}
\left(
q^n-mp^b(q^{n-1}-q^{n-2}+q^{n-3}) 
\right) 
\\&+&
\sum_{m=1}^{M_1}
\sum_{k=mk_0+1}^{mk_1}
\left(mq^{n-3}(q-p^b)+m-{k}\frac{q^{n-2}+1}{q+1}+\frac{q^n+1}{q+1}-1\right)
\\&+&
\sum_{m=M_1+1}^{M_2}
\sum_{k=mq+m-q}^{mk_1}
\left(mq^{n-3}(q-p^b)+m-{k}\frac{q^{n-2}+1}{q+1}+\frac{q^n+1}{q+1}-1\right) 
\\&+&
\sum_{m=M_2+1}^{M_3}
\sum_{ k=mq+m-q}^{q^2-q+m(q-p^b+1)}
\left(mq^{n-3}(q-p^b)+m-{k}\frac{q^{n-2}+1}{q+1}+\frac{q^n+1}{q+1}-1\right) 
\\
&=&M_1q^n-p^b(q^{n-1}-q^{n-2}+q^{n-3})M_1(M_1+1)/2 
  \\&+&
  \left(\frac{q^n+1}{q+1}-1\right) 
  \left(
  \sum_{m=1}^{M_1}
\sum_{k=mk_0+1}^{mk_1} 1 
+
\sum_{m=M_1+1}^{M_2}
\sum_{k=mq+m-q}^{mk_1}1 
+
\sum_{m=M_2+1}^{M_3}
\sum_{k=mq+m-q}^{q^2-q+m(q-p^b+1)}1 
\right) 
 \\&+&
(q^{n-3}(q-p^b)+1)
\left(
\sum_{m=1}^{M_1} m 
\sum_{k=mk_0+1}^{mk_1} 1 
+
\sum_{m=M_1+1}^{M_2} m 
\sum_{k=mq+m-q}^{mk_1}1 
+
\sum_{m=M_2+1}^{M_3} m 
\sum_{k=mq+m-q}^{q^2-q+m(q-p^b+1)}1 
\right) 
 \\&-&
\frac{q^{n-2}+1}{q+1}
\left(
\sum_{m=1}^{M_1}
\sum_{k=mk_0+1}^{mk_1} k 
+
\sum_{m=M_1+1}^{M_2}
\sum_{k=mq+m-q}^{mk_1}k 
+
\sum_{m=M_2+1}^{M_3}
\sum_{k=mq+m-q}^{q^2-q+m(q-p^b+1)}k 
\right) 
%
%
\\&=&M_1q^n-p^b(q^{n-1}-q^{n-2}+q^{n-3})M_1(M_1+1)/2 
  \\&+&
  \left(\frac{q^n+1}{q+1}-1\right) 
  \left(
A
+
B
+
C
\right) 
 \\&+&
(q^{n-3}(q-p^b) +1)
\left(
D
+
E
+
F
\right) 
 \\&-&
\frac{q^{n-2}+1}{q+1}
\left(
G
+
H
+
I
\right),
\end{eqnarray*}
where
{\footnotesize
\begin{eqnarray*}
A&=&\sum_{m=1}^{M_1}\sum_{k=mk_0+1}^{mk_1}1=
\sum_{m=1}^{M_1}(mk_1-mk_0)=
(k_1-k_0)\frac{M_1(M_1+1)}{2}\\
B&=&\sum_{m=M_1+1}^{M_2}\sum_{k=mq+m-q}^{mk_1}1=
\sum_{m=M_1+1}^{M_2}(mk_1-mq-m+q+1)\\&=&
(q+1)(M_2-M_1)+(k_1-q-1)\sum_{m=M_1+1}^{M_2}m\\&=&
(q+1)(M_2-M_1)+(k_1-q-1)(\frac{M_2(M_2+1)-M_1(M_1+1)}{2})\\
C&=&\sum_{m=M_2+1}^{M_3}(q^2-mp^b+1)=
(q^2+1)(M_3-M_2)-p^b(\frac{M_3(M_3+1)-M_2(M_2+1)}{2})\\
D&=&\sum_{m=1}^{M_1}m(mk_1-mk_0)=
(k_1-k_0)\sum_{m=1}^{M_1}m^2=
(k_1-k_0)(\frac{M_1(M_1+1)(2M_1+1)}{6})\\
E&=&\sum_{m=M_1+1}^{M_2}m(mk_1-mq-m+q+1)=
(q+1)\sum_{m=M_1+1}^{M_2}m+(k_1-q-1)\sum_{m=M_1+1}^{M_2}m^2\\&=&
(q+1)(\frac{M_2(M_2+1)-M_1(M_1+1)}{2})+(k_1-q-1)\frac{M_2(M_2+1)(2M_2+1)-M_1(M_1+1)(2M_1+1)}{6}\\
F&=&\sum_{m=M_2+1}^{M_3}m(q^2-mp^b+1)\\&=&
(q^2+1)(\frac{M_3(M_3+1)-M_2(M_2+1)}{2})-p^b(\frac{M_3(M_3+1)(2M_3+1)-M_2(M_2+1)(2M_2+1)}{6})\\
G&=&\sum_{m=1}^{M_1}(\frac{mk_1(mk_1+1)-mk_0(mk_0+1)}{2})=
\frac{k_1-k_0}{2}\sum_{m=1}^{M_1}m+\frac{k_1^2-k_0^2}{2}\sum_{m=1}^{M_1}m^2\\&=&
\frac{k_1-k_0}{2}\frac{M_1(M_1+1)}{2}+\frac{k_1^2-k_0^2}{2}\frac{M_1(M_1+1)(2M_1+1)}{6}\\&=&
\frac{(k_1-k_0)m(m+1)}{4}+\frac{(k_1^2-k_0^2)M_1(M_1+1)(2M_1+1)}{12}\\
H&=&\sum_{m=M_1+1}^{M_2}(\frac{mk_1(mk_1+1)-(m(q+1)-q-1)(m(q+1)-q)}{2})\\&=&
 -\frac{q(q+1)}{2}(M_2-M_1)+(\frac{k_1+(2q+1)(q+1)}{2})\frac{M_2(M_2+1)-M_1(M_1+1)}{2}\\&&+\frac{k_1^2-(q+1)^2}{2}\frac{M_2(M_2+1)(2M_2+1)-M_1(M_1+1)(2M_1+1)}{6}\\&=&
-\frac{q(q+1)}{2}(M_2-M_1)+\frac{(k_1+(2q+1)(q+1))(M_2(M_2+1)-M_1(M_1+1))}{4}\\&&+\frac{(k_1^2-(q+1)^2)(M_2(M_2+1)(2M_2+1)-M_1(M_1+1)(2M_1+1))}{12}
\\
I&=&\sum_{m=M_2+1}^{M_3}\sum_{k=mq+m-q}^{q^2-q+m(q-p^b+1)}k\\&=&
\sum_{m=M_2+1}^{M_3}\left(\frac{(q^2-q+m(q-p^b+1))(q^2-q+1+m(q-p^b+1))-(m(q+1)-q-1)(m(q+1)-q)}{2}\right)\\&=&
\frac{(q^2-q)(q^2-q+1)-q(q+1)}{2}(M_3-M_2)\\&&+
\frac{((q-p^b+1)(2q^2-2q+1)+(q+1)(2q+1))
(M_3(M_3+1)-M_2(M_2+1))}{4}\\&&
\frac{((q-p^b+1)^2-(q+1)^2)(M_3(M_3+1)(2M_3+1)-M_2(M_2+1)(2M_2+1))}{12}
\end{eqnarray*}
}
A {\sc Sage} simplification of this leads to 
\begin{eqnarray*}
&=&
  1/2 q^{n+1}p^{a-b} - 1/2q^n - 1/2q^2p^{a-b} + 1/2qp^{a-b}
  =
  \frac{q^{n+1}p^{a-b} - q^n - q^2p^{a-b} + qp^{a-b}}{2}
\end{eqnarray*}
\end{proof}

We remark here that the result does not vary if we
replace 
$M_2=(q-1)p^{a-b}$ by $M_2'=(q-1)p^{a-b}-1$.

\subsection*{Acknowledgments}

The first author was supported by Fundação de amparo a pesquisa de Minas Gerais under grant FAPEMIG APQ-00696-18.
The second author was supported by the Spanish government under grant TIN2016-80250-R and by the Catalan government under grant 2017 SGR 00705.

\end{document}